\documentclass{amsart}
\usepackage{graphicx}
\usepackage[all]{xy}
\usepackage{amscd}
\usepackage{amssymb}
\usepackage{amsmath}
\usepackage{amsthm}
\usepackage{amsfonts}
\usepackage{graphics}
\usepackage{epsfig,psfrag}
\usepackage[english]{babel}
\usepackage{latexsym}
\usepackage{mathrsfs}

\theoremstyle{definition}

\theoremstyle{remark}

\begin{document}

\title[A remark on the counterexample to the unknotting number conjecture]
{A remark on the counterexample to the unknotting number conjecture}

\author{Chao Wang}
\address{School of Mathematical Sciences, Key Laboratory of MEA(Ministry of Education) \& Shanghai Key Laboratory of PMMP, East China Normal University, Shanghai 200241, China}
\email{chao\_{}wang\_{}1987@126.com}

\author{Yimu Zhang}
\email{zym534685421@126.com}



\thanks{The first author is supported by National Natural Science Foundation of China (NSFC), grant Nos. 12131009 and 12371067, and Science and Technology Commission of Shanghai Municipality (STCSM), grant No. 22DZ2229014.}

\begin{abstract}
By using Snappy, M. Brittenham and S. Hermiller discovered a very surprising example that $u(7_1\#\overline{7_1})\leq 5<6=u(7_1)+u(\overline{7_1})$, where $7_1$ is the $(2,7)$-torus knot and $\overline{7_1}$ is its mirror image. Based on their work, we give a direct verification of this fact.
\end{abstract}

\date{}
\maketitle

\begin{figure}[h]
\includegraphics{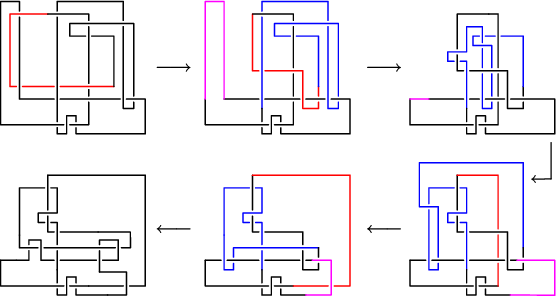}
\caption{Deformations between two diagrams of $K15n81556$.}\label{fig:K15}
\end{figure}
\begin{figure}[h]
\includegraphics{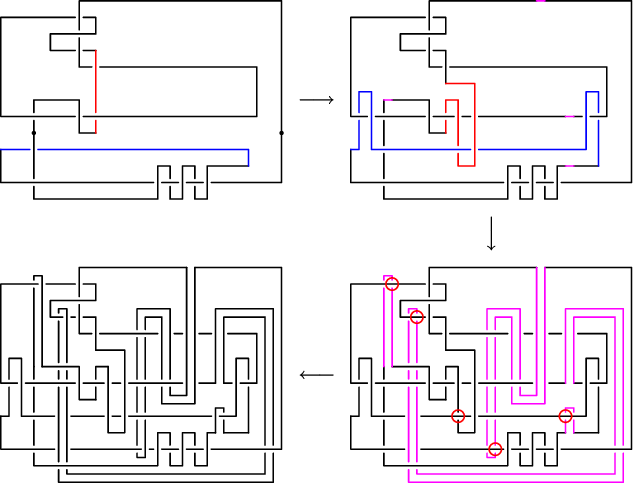}
\caption{Transformations from $7_1\#\overline{7_1}$ to the unknot.}\label{fig:Snake}
\end{figure}

The unknotting number $u(K)$ of a knot $K$ is the minimal number of crossing changes so that one can transform $K$ to the unknot. It is a very old conjecture that $u(K)$ is additive under the connected sum operation of knots. This can also be found as question 1.69(B) in Kirby's problem list \cite{K}. However, by using Snappy, M. Brittenham and S. Hermiller discovered a very surprising example that
\[u(7_1\#\overline{7_1})\leq 5<6=u(7_1)+u(\overline{7_1}),\]
where $7_1$ is the $(2,7)$-torus knot and $\overline{7_1}$ is its mirror image. See \cite{BH}. A key point in their proof of this fact is that they found $7_1\#\overline{7_1}$ can be transformed by $3$ crossing changes to the knot $K15n81556$, which can be transformed by $2$ crossing changes to the unknot. We noticed that the two diagrams of $K15n81556$ in \cite{BH} actually do not represent the same knot, but a chiral knot and its mirror image. This can be verified by the Jones polynomial and Figure~\ref{fig:K15}. Based on this fact, we can provide a direct verification of the inequality $u(7_1\#\overline{7_1})\leq 5$, as shown in Figure~\ref{fig:Snake}.

Note that the Figures in \cite{BH} together with Figure~\ref{fig:K15} can already provide a direct verification. By reversing this process, now we give a specific diagram of $7_1\#\overline{7_1}$ so that the five crossing changes can be done simultaneously. Of course, we will need a little more work to verify that the final result is indeed the unknot. This can be done by a sequence of Reidemeister moves, or by various other tools. We prefer to leave it to the readers as an interesting game.


\bibliographystyle{amsalpha}

\end{document}